\documentclass[11pt]{article}
\usepackage{amssymb}
\usepackage{graphicx}
\usepackage{amsmath} 
\usepackage{float}
\usepackage{multirow}
\newtheorem{theorem}{Theorem}

\newtheorem{definition}[theorem]{Definition}

\newtheorem{lemma}[theorem]{Lemma}
\newtheorem{proposition}[theorem]{Proposition}
\newtheorem{remark}[theorem]{Remark}
\newcommand{\R}{\mathbb{R}}
\newcommand{\He}{\mathrm{Hess}}
\newenvironment{proof}[1][Proof]{\textbf{#1.} }{\ \rule{0.5em}{0.5em}}
\date{}
\begin{document}
\title{{\bf On the realization problem of plane real algebraic curves as Hessian curves}}
\author{
Angelito Camacho Calder\'on\thanks{
Work supported by CONACyT } 
 $ $ and
Adriana Ortiz
Rodr\'{\i}guez\thanks{
Work partially supported by DGAPA-UNAM grant PAPIIT-IN108112 and N103010} }
\maketitle
\begin{abstract}
\noindent The Hessian Topology is a subject having interesting 
relations with several areas, for instance, differential geometry,
implicit differential equations, analysis and singularity theory. 
In this article we study the problem of realization of a
real plane curve as the Hessian curve of a smooth function.
The plane curves we consider are constituted either by only outer 
ovals or inner ovals. We prove that some of such curves are realizable 
as Hessian curves.
\end{abstract}


\noindent {\small {\it MS classification: 53A15; 53A05.}}

\section*{Introduction}
It is well known that the points of a generic smooth surface
in $\mathbb R^3$ ($\mathbb R$P$^3$) are classified as elliptic, 
parabolic and hyperbolic.
The geometric structure of such a surface is described by
the behavior of the sets conformed by these points, that is,
the set of parabolic points is a smooth curve called 
{\it parabolic curve} whose complement is a disjoint union 
of elliptic and hyperbolic domains.
A problem that has been studied for a long time is the description 
of this geometric structure, in particular the topology
of the parabolic curve \cite{salmon}, \cite{ker}, \cite{banch}, 
\cite{arn1}, \cite{arn2}, 
among others. When the surface is expressed locally as the graph $\,
z=f(x,y)\,$ of a smooth function $f$, the parabolic curve is the image
under $f$, of the {\it Hessian curve of $f$}, which is the set of zeros of the
Hessian polynomial of $f$, Hess $f:=f_{xx}f_{yy}-f_{xy}^2$. Moreover, if the surface
is the graph of a smooth function defined on the whole plane, then this 
local description is global, and so, to understand the behavior of the parabolic
curve in this case, it is enough to study the Hessian curve.

Some topics concerning the Hessian Topology \cite{geoast}, \cite{arn3} 
(problems 2000-1, 2000-2, 2001-1, 2002-1) are based in the study of the 
following problems.
How many compact connected
components can have the (smooth) parabolic curve of the graph of a real polynomial
of degree $n$ in two variables? and, how are they distributed?
\cite{or}, \cite{orso}, \cite{bb}, \cite{mos1}, \cite{mos2}.
Another such a problem is to determine the number of connected components 
of the parabolic curve of a smooth algebraic surface of fixed degree $n$ 
in the real projective space \cite{segre}, \cite{or}, \cite{bb}.
The realization problem of a smooth plane curve as the Hessian curve of
a smooth function is another such a problem of this subject, which is treated here. 
The meaning of all these problems is to find out if the property of being
Hessian curve imposes conditions on the topology and geometry of curves
and surfaces.   

A connected component of a compact smooth plane curve is called
{\it outer oval} if it is contained in the outside component of 
the complement of each of the connected components of the curve.
In this paper we give a family of real polynomials 
$\{f_k\}_{k\geq 2}\,$ 
such that the Hessian curve of $f_k$ consists exactly of $k-2$ outer ovals.
Jointly, in Theorem \ref{teo1}, we give a detailed description of the 
geometric structure of the graph of each one of these polynomials.

In Theorem \ref{teo-par}, for each even natural number $k$, 
we present real polynomials of degree $\,k\,$ such that its Hessian curve is 
formed exactly by $k-2$ concentric circles. We also describe
the geometric structure of the graph of these polynomials.
In an analogue form, we contribute to the odd case with a family of real functions
$\{f_r\}_{r\geq 1}\,$ such that the Hessian curve of $f_r$ 
are constituted by $2r-1$ concentric circles, Theorem \ref{teo-impar}.

To conclude, we show in Proposition \ref{hessianas-afines}, 
that if a smooth plane curve $\,C\,$ is 
equivalent by an affine transformation on the plane to any of the
previous Hessian curves, then $C$ is also a Hessian curve.

\section*{Preliminaries}

Points on a smooth surface (not necessarily generic) in $\mathbb R^3$ 
($\mathbb R$P$^3$) can be classified 
in terms of the maximal order of contact
of the tangent lines to the surface at each point.
A point $p$ of a surface
is {\it elliptic} if all lines
tangent to the surface at $p$ have a contact of order
$2$ with the surface at that point. 
The tangent lines having an order of contact at least three with the
surface are known as {\it asymptotic lines}.
A point $p$ is {\it hyperbolic} if it has exactly two asymptotic lines. 
We say that an hyperbolic point is an {\it inflexion point} if it has at least
one asymptotic line with order of contact greater than $3$.

The points on the surface that are not elliptic or 
hyperbolic will be called {\it parabolic points} and the
set constituted by them will be referred to as {\it parabolic curve}.
This curve is formed by the following types of points (some of them may 
not appear).
A {\it generic parabolic point} is a point which has exactly one 
asymptotic line and the order of contact of this line with the surface is $3$.
A {\it special parabolic point} (other authors call it Gaussian cusp or
godron) is a point at which the parabolic curve is smooth locally, the
4-jet (without constant neither linear terms) of the function defining 
the surface at this point is not a perfect square and moreover, 
it has exactly one asymptotic line whose order of contact is at least 4.
Finally, a {\it degenerate parabolic point} is a parabolic point 
that is not generic or special. It is important to mention that all
these type of points are invariant under the action of the affine group
on  3-space.

The directions determined on the $xy$-plane of
the asymptotic lines under the projection $\pi: \mathbb R^3 \rightarrow 
\mathbb R^2,\ (x,y,z) \mapsto (x,y)$,
are the solutions of the quadratic differential equation:
$$f_{xx}(x,y)dx^2+ 2f_{xy}(x,y)dxdy +f_{yy}(x,y)dy^2=0,$$
\noindent where the quadratic differential form on the left 
will be referred to as the {\it second fundamental form of $f$}. 
Its discriminant defined as 
$$\Delta_{II_f}= f_{xy}^2- f_{xx}f_{yy},$$ 

\noindent allows us to characterize the type of point in the graph of $f$. 
That is, $(p,f(p))$ is elliptic, parabolic or hyperbolic if  $\Delta_{II_f}(p)$ 
is negative, zero or positive, respectively.

\section*{Results and Proofs}
\section{Outer ovals}
We recall that a connected component of a compact smooth plane curve is called
{\it outer oval} if it is contained in the outside component of
the complement of each of the connected components of the curve.

\begin{definition}
{\rm A set of polynomials $\, l_1, \ldots ,l_n \in \R[x,y]\,$ of degree one
are in good position  if for each
$\,i=1,\ldots ,n,\,$ the straight line $\,l_i(x,y)$= $0$ contains no
critical point of the function $\,\prod_{j\neq i} l_j$.}
\end{definition}

\begin{theorem}\label{teo1}
Let $\,m,n\,$ two natural numbers.
Consider the polynomial
$$\, f(x,y) = \prod_{i=0}^m l_i(x,y) \prod_{j=0}^n {g_j}(x,y)$$
of degree $\,m+n+2\,$ given as the product
of polynomials of degree one, $ l_0(x,y)=y-ax,\, {g_0}(x,y)=y-bx,
\, l_i(x,y)=x-a_i, \, {g_j}(x,y)=x-b_j$, with $\,a, b, a_i, b_j \in \R,
\,i=1,\ldots, m, j=1,\ldots ,n,\, a\neq b\,$ and $\, a_m <
\cdots a_1<0<b_1<\cdots <b_n$.
If the set of polynomials $\, \{l_i(x,y), \, g_j(x,y):i=0,\ldots, m, j=0,\ldots ,n\}\, $
are in good position, then the Hessian curve of $f$ consists exactly of $\, m+n\,$
outer ovals. Moreover, the graph of $f$ has $\,3(m+n)\,$ special parabolic points and
its inflexion curve is the set $\,\{f(x,y)=0\}$.
\end{theorem}

\noindent\begin{proof}
In order to prove that the inflexion curve is the set $\,\{f(x,y)=0\}$,
we state the following Remark which is a consequence of the definition of $f(x,y)$.
\begin{remark}\label{lem2}
{\rm The restriction of $ f$ to any straight line
$m\in \R^2$ is a one variable real polynomial such that all its inflexion points are non degenerated.}
\end{remark}

\begin{lemma}\label{lem1}
The inflexion curve of the graph of $f$ is the set $\,\{f(x,y)=0\}$.
\end{lemma}

\noindent\begin{proof}
Because the order of contact of each straight line $\,l_i(x,y)= 0,$
$\,i=1,\ldots ,n,\,$ with the graph of $f$ is infinite then these lines are
contained in the inflexion curve. Inversely, suppose that $p$ is a point of
the inflexion curve and does not belong to the set $\,L:=\{f(x,y)=0\}$.
So, at least one asymptotic direction of $p$ has order of contact equal to
or greater than 4  with the graph of $f$. This implies that the restriction
of $f$ to the line, which is the projection to the $xy$-plane of
this asymptotic line, has a degenerated inflexion point. This contradicts
Remark \ref{lem2}.
\end{proof}

\begin{lemma}\label{lem4}
A point is in the intersection of the parabolic curve with the inflexion curve if and only if this point is a special parabolic point or a degenerated parabolic point.
\end{lemma}

\noindent\begin{proof}
If a point belongs to both curves, parabolic and inflexion,
it has an asymptotic line with order of contact at least 4 with the surface.
Then, the point is special parabolic or degenerated. Inversely, because a special
parabolic point satisfies that its asymptotic line has order of contact at least 4
with the surface, it lies in the inflexion curve. It only remains to prove that
any degenerated point lies in the inflexion curve. Note that a degenerated point
has one or an infinite number of asymptotic lines. For the first case, if the Hessian
curve is singular at this point, the asymptotic line has order of contact at
least 4 with the surface, so it lies in the inflexion curve. Remaining in the
 first case,  suppose that the Hessian curve is smooth at the point and
the 4-jet (without constant neither linear terms)
of the function is a perfect square. This implies that the point lies in
the inflexion curve. For the case in which every tangent line at the point
is an asymptotic line, we have that the second fundamental form of $f$ at the point
is identically zero, that is, each tangent line has order of contact at least 4
with the surface.
\end{proof}

\begin{lemma}\label{lem3}
The intersection of the parabolic curve with the inflexion curve
consists of $\,3(m+n)$ points.
\end{lemma}

\noindent\begin{proof}
After some straightforward calculus we have that the restriction of the
Hessian of $f$
to the straight line $\,x= a_i,\,$ for $i=1, ,\ldots ,m,\,$ is
\begin{eqnarray*}
\He f_{|_{x=a_i}}(x,y)&=& -(2y-(a+b)a_i)^{2}\Big(\prod_{j=1}^{n}(a_i-b_{j})\sum_{k=1}^{m}\prod_{s\neq k}^{m}(a_i-a_{s})\Big)^2 \\
&=&-(2y-(a+b)a_i)^{2} \,C_i^2,
\end{eqnarray*}
where $\, C_i= \prod_{j=1}^{n}(a_i-b_{j}) \sum_{k=1}^{m}\prod_{s\neq k}^{m}(a_i-a_{s})\,$  is a nonzero real number for each $\,i=1,\ldots, m.$

Then $\He f_{|_{x=a_i}}(x,y) = 0\,$ if and only if $\,y=\frac{(a+b)a_i}{2},
\,i=1,\ldots,m.$ The points $\left(a_i,\frac{(a+b)a_i}{2}\right),$ $i=1,\ldots, m,$
are between the line $l_0(x,y)=0$ and the line $g_0(x,y)=0$ because
$\, b a_i < \frac{(a+b)a_i}{2} <a a_i$. Moreover, the other points on the
lines  $l_{i}=0,\, i=1,\ldots, m,$ are hyperbolic points.

By a similar analysis for the lines $\,x= b_j\,$ for $j=1, ,\ldots ,n$
we have that the points $\left(b_j,\frac{(a+b)b_j}{2}\right),$ $j=1,\ldots, n,$ are
parabolic and they are in the interior of the compact segment
$[a b_j,b b_j]$ of the line $\,x= b_j\,$.

Now, let us analyze the intersection of the parabolic curve with the straight lines $l_0(x,y)=0$ and $g_0(x,y)=0$.
Replacing $\,y=ax$ in the Hessian polynomial of $f$ we obtain
\begin{eqnarray*}
\He f (x,y)|_{y=ax}= -(a-b)^2\left[\prod_{i=1}^{m}(x-a_i)\prod_{j=1}^{n}(x-b_j)+ \right.\qquad\qquad\qquad\\
\left. + \,x\left(\sum_{l=1}^{m}\prod_{i\neq l}^{m}(x-a_i)\prod_{j=1}^{n}(x-b_j)+
\sum_{k=1}^{n}\prod_{j\neq k}^{n}(x-b_j)\prod_{i=1}^{m}(x-a_i)\right)\right]^2.
\end{eqnarray*}

Note that the expression being inside of the square brackets is the first
derivative of the one variable polynomial
$\displaystyle g(x)=x\prod_{i=1}^{m}(x-a_i)\prod_{j=1}^{n}(x-b_j)$ of degree
$m+n+1$ which has $m+n+1$ simple zeros and  $m+n$ simple critical points.
This implies that $\He f (x,y)|_{y=ax}$ has $m+n$ simple real roots.

By an analogous analysis for the straight line $\,y=bx$ we obtain that
$\He f (x,y)|_{y=bx}$ has $m+n$ simple real roots.
\end{proof}

\begin{lemma}\label{lem6}
If the set of polynomials $\, \{l_i(x,y), \, g_j(x,y):i=0,\ldots, m, j=0,\ldots ,n\}\,$
are in good position, then $\,\{f(x,y)=0\}$ contains no degenerated parabolic point.
\end{lemma}

Before starting with the proof of this Lemma we state the following Remark which is
obtained by a suitable grouping of terms of the Hessian polynomial of $f$.
\begin{remark}\label{agrupa-hess}
{\rm The Hessian polynomial of $f$ can be written as
$$\He f(x,y) =\beta(x)\left(y-\frac{(a+b)x}{2}\right)^2+\alpha(x),$$ where
{\small \begin{eqnarray*}
\alpha(x)=-(a-b)^2\prod_{i=1}^{m}(x-a_i)\prod_{j=1}^{n}(x-b_j)\left[ \prod_{i=1}^{m}(x-a_i)\prod_{j=1}^{n}(x-b_j)\right.\qquad\qquad\quad\\
+2x\left(\sum_{l=1}^{m}\prod_{i\neq l}^{m}(x-a_i)\prod_{j=1}^{n}(x-b_j)
+\sum_{k=1}^{n}\prod_{j\neq k}^{n}(x-b_j)\prod_{i=1}^{m}(x-a_i)\right)\qquad\quad\,\,\,\\
+\frac{x^2}{2}\left(\sum_{l=1}^{m}\sum_{s\neq l}^{m}\prod_{i\neq l,s}^{m}(x-a_i)\prod_{j=1}^{n}(x-b_j)+\sum_{k=1}^{n}\sum_{t\neq k}^{n}\prod_{j\neq k,t}^{n}(x-b_j)\prod_{i=1}^{m}(x-a_i)\right.\\
\left.\left.+2\sum_{k=1}^{n}\prod_{j\neq k}^{n}(x-b_j)\sum_{l=1}^{m}\prod_{i\neq l}^{m}(x-a_i)\right)\right],\hskip 5.3cm
\end{eqnarray*}}
{\small \begin{eqnarray*}
\beta(x)=2\left(y-\frac{(a+b)x}{2}\right)^2\left[\prod_{i=1}^{m}(x-a_i)\prod_{j=1}^{n}
(x-b_j)\bigg(\sum_{l=1}^{m}\sum_{s\neq l}^{m}\prod_{i\neq l,s}^{m}(x-a_i)\prod_{j=1}^{n}(x-b_j)\right.\\
+\sum_{k=1}^{n}\sum_{t\neq k}^{n}\prod_{j\neq k,t}^{n}(x-b_j)\prod_{i=1}^{m}(x-a_i)+2\sum_{k=1}^{n}\prod_{j\neq k}^{n}(x-b_j)\sum_{l=1}^{m}\prod_{i\neq l}^{m}(x-a_i)\bigg)\qquad\quad\\
\left.-2\left(\sum_{l=1}^{m}\prod_{i\neq l}^{m}(x-a_i)\prod_{j=1}^{n}(x-b_j)+
\sum_{k=1}^{n}\prod_{j\neq k}^{n}(x-b_j)\prod_{i=1}^{m}(x-a_i)\right)\right].\qquad\qquad
\end{eqnarray*}}
}
\end{remark}

\noindent\begin{proof}(Lemma \ref{lem6}).
We shall prove that the points lying in both curves, parabolic and inflexion,
are all special points.

The next calculus shows that the second fundamental form of $f$ at each point  $q\in \{f(x,y)=0\}\,$ is nonzero. Let suppose that $\, q\in \{l_{i}(x,y)=0\}$
for some $i\in \{0,\ldots,m\}$ (the following arguments are valid for the case that
$q$ belongs to some $l_{i}(x,y)=0$). So, $\,f(x,y)=l(x,y) l_{i}(x,y)$ where $ l(x,y)=\prod_{k=0,k\neq i}^{m}l_{k}(x,y) \prod_{j=0}^n g_{j}(x,y).$ Then,
\begin{eqnarray*}
f_{xx}(q)=2\frac{\partial l_{i}}{\partial x}(q) l_{x}(q),\qquad
f_{yy}(q)=2\frac{\partial l_{i}}{\partial y}(q) l_{y}(q),\\
f_{xy}(q)=\frac{\partial l_{i}}{\partial x}(q) l_{y}(q)+\frac{\partial l_{i}}{\partial y} l_{x}(q). \qquad\qquad
\end{eqnarray*}

The hypothesis of good position implies that $q$ can not be critical point of $l$.
So, $l_{x}(q)\neq 0$ or $l_{y}(q)\neq 0$ and the second fundamental form of
$f$ at $q$ is non degenerated. That is, each point has exactly one asymptotic
line.

Now, let $p$ be a point in the intersection of the parabolic curve with the
inflexion curve. We shall prove that the Hessian curve at $p$ is non singular.
The proof will be in two cases.
\begin{itemize}
\item  The point $p$ lies in any of the lines $\,l_{0}(x,y)=0\,$ or $\,g_{0}(x,y)=0$.
After an appropriate affine transformation of the $xy$-plane, we can suppose that
such a line is the line $y=0$. Denote by $\, x_1, \ldots ,x_{m+n+1}$
the intersection points of the line $y=0$ with all other lines. Since
$ f(x,y)=y \prod_{k=1}^{m+n+1} r_k(x,y)$, where $\,r_k(x,y)$ are real polynomials
of degree one, then,
$$f_y(x,y)|_{y=0} = \prod_{k=1}^{m+n+1} r_k(x,y)|_{y=0}.$$
Define $\,f_y(x,y)|_{y=0} = F(x)$. The polynomial $F(x)$ is of degree $m+n+1$ and
its roots are real and simple because they are the points $\, x_1, \ldots ,x_{m+n+1}$.
Moreover, it has $m+n$ non degenerated critical points, that is, the function
$\,F'(x) = f_{xy}(x,y)|_{y=0} \,$
has $m+n$ simple real zeroes: exactly one in the interior of each segment $\,[x_i,x_{i+1}],
\, i=1, \ldots ,m+n$ and $\,f_{xxy}(x,y)|_{y=0} \neq 0\,$ at these roots.
Note that
$$(\He f(x,y))|_{y=0} = - (f_{xy}(x,y))^2|_{y=0} = - (F'(x))^2.$$
So, the point $p$, after applying the affine transformation,
has the form $\tilde p = (\tilde x, 0)$, where $\tilde x$ is one of the critical points of $F$. Moreover, $\,f_{xy}(x,y),\,f_{xx}(x,y)$ are zero at this point.
Finally, since the second fundamental form of $f$ at $\tilde p$ is non zero and
$\,f_{xxy}(\tilde p)\neq 0$, then
\begin{eqnarray*}
\left(\frac{\partial}{\partial y} \He f\right)(\tilde p) = (f_{xx}f_{yyy}+f_{xxy}f_{yy}-2f_{xy}
f_{xyy})|_{\tilde p}\\
= f_{xxy}(\tilde p)f_{yy}(\tilde p) \neq 0.\hskip 2.39cm
\end{eqnarray*}

\item The point $p$ lies in one of the lines $\,l_{1}(x,y)=0,\ldots ,l_{m}(x,y)=0,$
$g_{1}(x,y)=0,\ldots ,g_{n}(x,y)=0.$
In the proof of Lemma \ref{lem3} it is shown that  $p$ has the form
$\left(x_0,\frac{(a+b)x_0}{2}\right)\,$  where $\,x_0\in\{a_1,\ldots ,a_m,b_1,\ldots ,b_n\}$.
So, by Remark \ref{agrupa-hess} we have
$\,\frac{\partial}{\partial y}\He f(p)= \alpha'(x_0) \neq 0$.
\end{itemize}

Note that the asymptotic line at each one of these points is a line of the set
$\, \{l_i(x,y)=0, \, g_j(x,y)=0 :i=0,\ldots, m, j=0,\ldots ,n\}$
because the parabolic curve is tangent to this set.

It only remains to prove that the 4-jet (without constant neither linear terms)
of the function at these points is free square.
For the points $\left(a_i,\frac{(a+b)a_i}{2}\right),$
this assertion follows by noting that the coefficient
of the $xy^2$ term is nonzero while the coefficient of $y^4$ is zero.
\end{proof}

\medskip
In order to prove the compacity of the parabolic curve consider the bounded connected components of the complement
of $\,\{f(x,y)=0\}\,$ in the $xy$-plane. Take its closure of each one
and the union of them. Let us denote by $W$ the complement of this union.
\begin{lemma}\label{lemas1}
The points belonging to the set $W$ are hyperbolic.
\end{lemma}

\noindent\begin{proof}
We partition the set $W$ in 4 regions which are
represented in Figure \ref{fig1}.

\begin{figure}[H]
\begin{center}
\includegraphics[width=5cm]{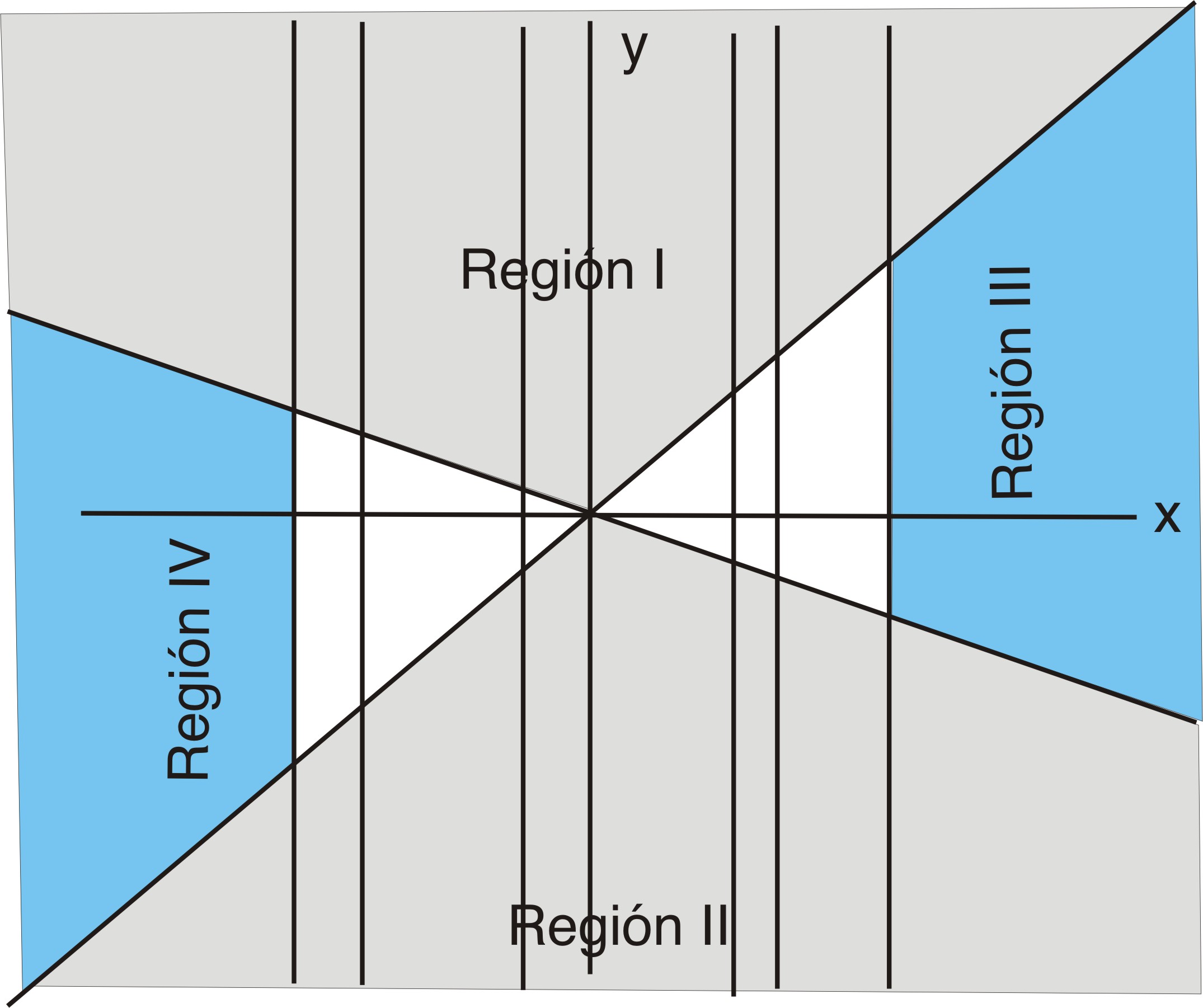}
\caption{Partition of the set $W$ in 4 regions}
\label{fig1}
\end{center}
\end{figure}

Consider Region I. It was shown in the proof of Lemma \ref{lem3}
that the points lying in some of the straight lines
$\,l_{i}(x,y)=0, \,{g_{j}}(x,y)=0$, with  $\,i=0,\ldots,m, \,j=0,\ldots,n\,$
of Region I (or II) are hyperbolic. So, it
only remains to show that the points of Region I and belonging to
the complement of $\,\{f(x,y)=0\}\,$ are hyperbolic.
Let us denote by $q$ such a point.
Consider the straight line $m$ on the $xy$-plane parallel to the $y$-axis
and passing through $q$. The restriction $\, f|_m\,$ of $\, f$ to this line
is a one variable polynomial of degree 2 with two real roots
(the intersections with the straight lines $\,l_{0}(x,y)=0, \,g_{0}(x,y)=0$).
So, $f(q)$ and the derivative of $\, f|_m\,$ at $q$ are positive or negative.
That is $\, f|_m\,$ is convex or concave at $q$, respectively.

Suppose that $f(q)$ is positive. Denote by $t$ the straight line tangent to the level curve $\,f(x,y) = f(q)$ at $q$. Because the derivative of
$\, f|_t\,$ is zero, the line $t$ is different from the line $m$.
The function $\, f|_t\,$ has a non degenerated critical
point at $q$ which is a maximum. So, the function $\, f|_t\,$ is concave at $q$. The two previous conclusions imply that $(q,f(q))$ is a hyperbolic point.
The case $\,f(q)<0\,$ is analogous.

The analysis for Region II is similar.

Now, we consider Region III which is the unbounded connected component
determined by the lines $\,x=b_n, \,y= ax\,$ and $\,y=bx$.
Let $q$ be a point in such a set. Consider the straight lines $\,m, l\,$
parallels to the lines $g_0(x,y)=0,\, l_0(x,y)=0$, respectively, and passing through $q$.
The lines $\,m,\, l$ divide the $xy$-plane in four sectors, one of which,
denoted by $A$, is totally contained in Region III.
Take a straight line $t$ such that it passes through $q$, it divides the
sector $A$ in two unbounded sectors and it is not parallel to any line of $L$.
So, the roots of $\,f|_t\,$ are in the left side of $q$.

If $\,f(q)>0$, then the derivative of $\,f|_t\,$ at $q$ is positive.
So $\,f|_t\,$ is convex at $q$.
On the other hand, consider the straight line passing through $q$ which is parallel
to the $y$-axis. Denote it by $r$. So, $\,f|_r\,$ is concave at $q$ because
it is a one variable polynomial of degree two and $\,f(q)>0$. Then, the point $(q,f(q))$ is hyperbolic.
The case for Region IV is similar.
\end{proof}

\medskip
Since there is one special parabolic point at each compact segment
of the straight lines $g_0(x,y)=0,\, l_0(x,y)=0$, there is at least
one connected component of the Hessian curve on each compact polygone
delimited by the lines $g_j(x,y)=0,\, l_i(x,y)=0, \, i=0,\ldots,m,\,
j=0,\ldots,n.$ Because there are $m+n$ such as compact polygons, then
there are at least $m+n$ connected components of the Hessian curve.
Now, the goal is to prove that the Hessian curve has exactly $m+n$
connected components. In fact, we shall show that the Hessian curve
has at most $2(m+n)$ vertical tangent lines since this implies that
the Hessian curve has at most $m+n$ connected components.

Let us prove it.
Consider the expression of the Hessian polynomial of $f$ given in Remark
\ref{agrupa-hess}. The parabolic points having a vertical tangent line to the
parabolic curve are given by
$$\Big\{ \frac{\partial}{\partial y} \mbox{Hess} f(x,y)=0\Big\}
\cap \Big\{\mbox{Hess} f(x,y)=0\Big\}.$$
Note that if $x_0$ is a real root of the polynomial $\beta(x)$ then
$\alpha(x_0)\neq 0$ because in another case the points of the line
$x=x_0$ would be parabolic and it contradicts Lemma \ref{lemas1}.
So, the parabolic points having a vertical tangent line are the points
$\,\left(\tilde{x_0},\frac{(a+b)\tilde{x_0}}{2}\right)\,$ such that $\tilde{x_0}$ is a real root
of $\alpha(x)$. The proof concludes by noting that
$\alpha(x)$ is a one variable real polynomial of degree $2(m+n)$.
\end{proof} 
\section{Concentric circles. Even case}
\begin{theorem}\label{teo-par}
Consider the real polynomial of degree $\,2n$
$$\,f(x,y)=\prod_{i=1}^{n}(x^2+y^2-m_{i}^2)\,$$ 
with $\,0<m_{1}<m_{2}<\cdots<m_n,\,$ $\,m_{i}\in \R, i=1,\ldots,n$.  
Then, its Hessian curve is non singular and it is formed by $2(n-1)$
concentric circles. Moreover, the unbounded connected component of the complement 
of the Hessian curve on the $xy$-plane is elliptic.
\end{theorem}

\noindent\begin{proof}
After a straightforward calculus the Hessian polynomial of $f$ is
$$\He (f)=4 \,s(x,y)\, t(x,y),\quad \mbox{where}$$
\begin{eqnarray*}
s(x,y)=\sum_{j=1}^{n}\prod_{i\neq j}^{n}(x^2+y^2-m_{i}^2)\,\quad \mbox{ and}\hskip 2.9cm\\
t(x,y)=\sum_{j=1}^{n}\prod_{i\neq j}^{n}(x^2+y^2-m_{i}^2)+2(x^2+y^2)\sum_{j=1}^{n}\sum_{l\neq j}^{n}\prod_{i\neq j,l}^{n}(x^2+y^2-m_{i}^2).
\end{eqnarray*}

Note that the graph of the Hessian 
of $f$, $\, z-\He f(x,y)=0,$ is the rotation about the $z$-axis of the graph of the function
$\,4\tilde{s}(x) \tilde{t}(x),\,$ where 
$$\widetilde{s}(x)=\sum_{j=1}^{n}\prod_{i\neq j}^{n}(x^2-m_{i}^2)\, \mbox{ and}\qquad$$
$$\widetilde{t}(x)=\sum_{j=1}^{n}\prod_{i\neq j}^{n}(x^2-m_{i}^2)+2(x^2)\sum_{j=1}^{n}\sum_{l\neq j}^{n}\prod_{i\neq j,l}^{n}(x^2-m_{i}^2).$$
So, the solution curves of $\, s(x,y)=0$ $(\,t(x,y)=0)\,$ correspond to the non negative
real roots of $\widetilde{s}(x)$ $(\,\widetilde{t}(x))$. We shall prove now that 
$\widetilde{s}(x)$ and $\widetilde{t}(x)$ have $n-1$ simple positive real roots, each one. 
Define the one variable real polynomial
$$\widetilde{f}(x):=\prod_{i=1}^{n}(x^2-m_{i}^2).$$ 
The degree of this polynomial is $\,2n$, it has $2n$ non zero simple roots,
$2n-1$ non degenerated critical points and $2n-2$ simple inflexion points.
Because 
$$\widetilde{f}_{x}(x)=2x\sum_{j=1}^{n}\prod_{i\neq j}^{n}(x^2-m_{i}^2)=2x\widetilde{s}(x),$$
then $\,\widetilde{s}(x)$ has $2n-2$ non zero simple real roots from which, $n-1$
are positive by symmetry. Moreover, since
$$\widetilde{f}_{xx}(x)=2\sum_{j=1}^{n}\prod_{i\neq j}^{n}(x^2-m_{i}^2)+4x^2\sum_{j=1}^{n}\sum_{l\neq j}^{n}\prod_{i\neq j,l}^{n}(x^2-m_{i}^2)=\widetilde{t}(x),$$
then $\widetilde{t}(x)$ has $2n-2$ simple real roots (different from the roots of 
$\widetilde{s}(x)$). Because $\widetilde{f}$ is symmetric with respect to the origin, then $\widetilde{t}(x)$ has $n-1$ positive roots. Moreover, the unbounded connected component
of the complement of the Hessian curve on the $xy$-plane is elliptic because the 
functions $\tilde{s}$ and $\tilde{t}$ are positive if $\,x > m_n.$
\end{proof}

\section{Concentric circles. Odd case}
\begin{theorem} \label{teo-impar}
Consider the real valued differentiable function 
$$\,f(x,y)=\prod_{k=1}^{n}\frac{x^2+y^2-k^2}{x^2+y^2+1}.$$ 
Its Hessian curve is non singular and it is formed by $\,2n-1\,$ concentric
circles. Moreover, the unbounded connected component of the complement 
of the Hessian curve on the $xy$-plane is hyperbolic.
\end{theorem}

\noindent\begin{proof}
After a straightforward computation the Hessian polynomial of $f$ is
$$\He f(x,y)= \frac{4 \, s(x,y)\, t(x,y)}{(x^2+y^2+1)^{2n+3}},\,\, \mbox{ where}$$
\begin{eqnarray*}
s(x,y):= \sum_{j=1}^{n}\bigg((j^2+1)\prod_{i\neq j}^{n}(x^2+y^2-i^2)\bigg),\hskip 3.9cm\\
t(x,y):=(1-3x^2-3y^2)\sum_{j=1}^{n}\bigg((j^2+1)\prod_{i\neq j}^{n}(x^2+y^2-i^2)\bigg)+\hskip 1.0cm\\
\qquad\qquad +2(x^2+y^2)\sum_{j=1}^{n} \bigg( (j^2+1)\sum_{l\neq j}^{n}
\bigg( (l^2+1)\prod_{k\neq j,l}^{n}(x^2+y^2-k^2) \bigg) \bigg). 
\end{eqnarray*}
Note that the graph of the Hessian polynomial
of $f$, is the rotation about the $z$-axis of the graph of the function
$\,\frac{4\tilde{s}(x)\tilde{t}(x)}{(x^2+1)^{2n+3}}$, where
\begin{eqnarray*}
\tilde{s}(x)= \sum_{j=1}^{n}\bigg((j^2+1)\prod_{k\neq j}^{n}(x^2-k^2)\bigg),\,\,\quad\quad
\tilde{t}(x)=(1-3x^2)\sum_{j=1}^{n}\bigg((j^2+ \\
+1)\prod_{k\neq j}^n(x^2-k^2)\bigg)+2x^2\sum_{j=1}^{n}\bigg((j^2+1)\sum_{l\neq j}^{n}\bigg( (l^2+1)\prod_{k\neq j,l}^{n}(x^2-k^2)\bigg) \bigg).
\end{eqnarray*}
We analyze now the zeros of $\,\tilde{s}(x)\,$ and $\,\tilde{t}(x)$.
Consider the differentiable function $\,\widetilde{f}:\R \rightarrow \R$
$$\widetilde{f}(x)=\prod_{k=1}^{n}\frac{x^2-k^2}{x^2+1}.$$
This function has exactly $2n$ simple zeros from which $n$ are positive.
Moreover,
\begin{eqnarray*}
\widetilde{f}_{x}(x)=\frac{2x}{(x^2+1)^{n+1}}\sum_{j=1}^{n}\bigg((j^2+1)\prod_{i\neq j}^{n}(x^2-i^2)\bigg) = \frac{2x\,\widetilde{s}(x)}{(x^2+1)^{n+1}},\quad\qquad\\
\widetilde{f}_{xx}(x)=\frac{2}{(x^2+1)^{n+2}}\Bigg[(1-3x^2)\sum_{j=1}^{n}
\bigg((j^2+1)\prod_{i\neq j}^{n}(x^2-i^2)\bigg)+ \qquad\qquad\\
+ 2x^2\sum_{j=1}^{n}\bigg((j^2+1)\sum_{l\neq j}^{n}\bigg( (l^2+1)\prod_{k\neq j,l}^{n}(x^2-k^2)\bigg)\bigg)\Bigg]=
\frac{2\, \widetilde{t}(x)}{(x^2+1)^{n+2}}.
\end{eqnarray*}
Since $\widetilde{f}(x)$ has at least $2n-1$ critical points and the degree of $\tilde{s}(x)$ is $2n-2$ then, the polynomial $\tilde{s}(x)$ has exactly $2n-2$ simple
roots different from zero ($n-1$ are positive). So, the set $\,s(x,y)=0\,$ is a disjoint
union of $n-1$ concentric circles.

Because the critical points of $\widetilde{f}(x)$ are non degenerated and it is symmetric respect to the origin, then $\widetilde{f}(x)$ has at least $n-1$ positive inflexion points in the interval $(-n,n),$ exactly one point between two consecutive critical points.
That is, $\widetilde{t}(x)$ has $n-1$ positive real roots in the interval $[0,n)$.

Denote by $\,x_{0}$ the critical point of $\,\tilde{f}(x)$ lying in the interval $(n-1,n).$ Note that the sign of $\,\tilde{s}(n-1)= \left((n-1)^2+1\right)\prod_{k\neq n-1}^{n}\left((n-1)^2-k^2\right)\,$ is negative while the sign of
$\,\tilde{s}(n)= (n^2+1)\prod_{i=1}^{n-1}(n^2-i^2)\,$ is positive.
From these remarks and because the function $\,\widetilde{f}$ has exactly one 
critical point at each interval $\,(i,i+1),\, i=1,\ldots,n-1$, we have that $\,x_{0}$ is a minimum. Moreover, $\tilde{t}(x_{0})$ is positive.
After a straightforward calculation, it can be seen that $\,\tilde{t}(n^2)\,$ is negative. 
So, by the main value Theorem, there exists a real number $\,c\in (x_{0},n^2)\,$ 
such that $\tilde{t}(c)=0.$ It implies that $\tilde{t}(x)$ has $n$ positive real roots.
That is, the set $\,t(x,y)=0\,$ is a disjoint union of $n$ concentric circles different
from the circles described by the equation $\,s(x,y)=0.$
We conclude that the Hessian curve of $\,f\,$ consists of $2n-1$ concentric circles.
\end{proof}

\begin{proposition}\label{hessianas-afines}
If a smooth plane curve $\,C\,$ is equivalent
by an affine transformation of the plane, to a
Hessian curve, then $\,C\,$ is also a Hessian curve.
\end{proposition}

\noindent\begin{proof}
Suppose that the curve $C$ is defined by $\,g(x,y)=0,\,$
where $g$ is a smooth real valued function defined on the plane.
By hypothesis, the curve $C$ is equivalent by an affine
transformation to a Hessian curve, that is, $\,g(x,y)= 
(\mbox{Hess }f)\circ T(x,y),$ where $T:\R^2 \rightarrow \R^2$
is an affine transformation of the plane whose determinant 
(of the linear part) is denoted by $J$.
After a straightforward calculation, we have that 
$$\mbox{Hess }\left(\frac{f\circ T(x,y)}{J}\right) = (\mbox{Hess }f)\circ T(x,y).$$
So, the curve $C$ is also a Hessian curve.
\end{proof}

\noindent Adriana Ortiz Rodr\'{\i}guez, Instituto de Matem\'aticas,
Universidad Nacional Aut\'onoma de M\'exico,
Area de la Inv. Cient., Circuito Exterior,
C. U., M\'exico D.F 04510, M\'exico.\newline
e-mail: aortiz@matem.unam.mx

\noindent Angelito Camacho Calder\'on, Instituto de Matem\'aticas,
unidad Cuernavaca, Universidad Nacional Aut\'onoma de M\'exico,
Av. Universidad s/n. Col. Lomas de Chamilpa, c.p.62210, Cuernavaca, Morelos.\newline
e-mail: camacho@matem.unam.mx

\end{document}